\documentstyle[11pt]{article}
\begin{document}
\author{{Shiqiu Zheng$^{1, 2}$\thanks{E-mail: shiqiumath@163.com(S. Zheng).}\ , \ \ Shoumei Li$^{1}$\thanks{E-mail: lisma@bjut.edu.cn(S.Li).}}
  \\
\small(1, College of Applied Sciences, Beijing University of Technology, Beijing 100124, China)\\
\small(2, College of Sciences, North China University of Science and Technology, Tangshan 063009, China)\\
}
\date{}

\title{\textbf{Representation theorems for generators of reflected BSDEs with continuous and linear growth generators}\thanks{This work is supported by the Science and Technology Program of Tangshan (No. 13130203z).}}\maketitle

\textbf{Abstract:}\quad  In this paper, we establish a local representation theorem for generators of reflected backward stochastic differential equations (RBSDEs), whose generators are continuous with linear growth. It generalizes some known representation theorems for generators of backward stochastic differential equations (BSDEs). As some applications, a general converse comparison theorem for RBSDEs is obtained and some properties of RBSDEs are discussed.
\\


\section{Introduction}
The theory of backward stochastic differential equations (BSDEs) has gone through rapid development in many different research areas in recent 20 years (see Peng [1]). One of the important results of BSDEs theory is the representation theorem of generator, which establishes a relation between generators and solutions of BSDEs in limit form, can be used to resolve many problems in BSDEs theory. Representation theorem of generator is firstly obtained by Briand et al. [2], then generalized to the case that the generator $g$ only satisfy Lipschitz condition and applied to study the converse comparison theorems of BSDE and uniqueness theorem, translation invariance, convexity, etc, for $g$-expectation by Jiang [3][4]. Since then, the representation theorem of BSDEs is further studied under some more general conditions. For example, Jia [5] obtains a representation theorem for BSDEs with continuous and linear growth generators. Fan and Jiang [6] obtains a result for BSDEs with continuous and linear growth generators in the space of processes. Recently, Fan et al. [7] obtains a result for BSDEs whose generators are monotonic and polynomial growth in $y$ and linear growth in 	
$z.$

El Karoui et al. [8] introduced the notion of reflected BSDEs (RBSDEs). A RBSDE is a BSDE with an additional
continuous, increasing process $K$ in this equation to keep the solution
above a given continuous process $L$, called obstacle. RBSDEs are widely applied to the pricing of American options, mixed control, partial differential equations, etc. To interpret the relation between the generators and solutions of RBSDEs, the following problem is natural: Can we establish a representation theorem for generators of RBSDEs?

A basic difficulty to solve this problem is that the solution of RBSDEs is restricted by its obstacle. In this paper, we solve this problem using a localization method. Exactly, we establish a local representation theorem for generators of RBSDEs with continuous and linear growth generators, which generalizes some representation theorems for BSDEs (see Jiang [3] and Jia [5], etc) to RBSDEs case. Compared with the BSDEs case, our representation theorem of RBSDEs contains the increasing process $K,$ and is obtained in local space.

Converse comparison theorem of RBSDEs is firstly studied by Li and Tang [9], then by Li and Gu [10],  when generators $g$ are continuous in $t$  and satisfy Lipschitz condition. Using the representation theorem obtained in this paper, we obtain a general converse comparison theorem of RBSDEs, whose generators are only continuous with linear-growth in $(y,z)$. With the help of our representation theorem, we also discuss some properties of RBSDEs.

This paper is organized as follows. In the next section, we will recall the definition of RBSDEs. In section 3, a local representation theorem for generator is established under continuous and linear growth condition. In section 4, some applications of representation theorem are given.
\section{Preliminaries}
Let $(\Omega ,\cal{F},\mathit{P})$ be a complete probability space
carrying a $d$-dimensional standard Brownian motion ${{(B_t)}_{t\geq
0}}$, starting from $B_0=0$. Let $({\cal{F}}_t)_{t\geq 0}$ denote
the natural filtration generated by ${{(B_t)}_{t\geq 0}}$, augmented
by the $\mathit{P}$-null sets of ${\cal{F}}$, let $|z|$ denote its
Euclidean norm, for $\mathit{z}\in {\mathbf{R}}^n$,  let $T>0$ be a given real number. We define the following usual
spaces:

$L^p({\mathcal {F}}_t)=\{\xi:\ {\cal{F}}_t$-measurable
random variable; $E[|\xi|^p]<\infty\},\ t\in[0,T],\ p\geq1; $

${\mathcal{S}}^{\mathrm{2}}
 \left(0,{\mathit{T}};\mathbf{R}\right)=\{\psi:$ continuous predictable
process; $\|\psi\|_{\cal{S}^{\mathrm{2}}}^2=E\left
[\sup_{t\in[0,T]} |\psi _t|^2\right] <\infty \};$

${\cal{H}}^2(0,T;{\mathbf{R}}^d)=\{\psi:$ predictable
process;
$\|\psi\|_{{\cal{H}}^2}^2={E}\left[\int_0^T|\psi_t|^2dt\right]
<\infty \}.$

Let us consider a function $g$
$${g}\left( \omega ,t,y,z\right) : \Omega \times [0,T]\times \mathbf{%
R\times R}^{\mathit{d}}\longmapsto \mathbf{R}$$ such that
$\left(g(t,y,z)\right)_{t\in [0,T]}$ is progressively measurable for
each $(y,z)\in\mathbf{%
R\times R}^{\mathit{d}}$. In this paper, we
will make the following assumptions on $g$:

(A1). (Linear growth) There exists a constant
$\lambda>0$, and non-negative stochastic process $\gamma_t\in{\cal{H}^{\mathrm{2}}}(0,T;{\mathbf{R}})$ such that $P-a.s.,$ for each $t\in[0,T]$ and $(y,z)\in {\mathbf{ R\times R}}^{\mathit{d}},$
$$|{g}\left( t,y,z\right)|\leq \lambda\left(\gamma_t+|y|+|z|\right).$$

(A2). (Continuity) $P-a.s.,$ for each $t\in[0,T],$ $(y,z)\longmapsto g(t,y,z)$ is continuous.

(A3). $P$-$a.s.$, for each $(t,y)\in[0,T]\times{\mathbf{R}}^{\mathit{d}},$ $g(t,y,0)=0.$

\textbf{Definition 2.1} A RBSDE is associated with a terminal
condition $\xi\in L^2({\mathcal {F}}_T)$, a generator $g$,
a lower obstacle $\{L_t\}_{0\leq t\leq T},$ which is a continuous progressively measurable real-valued process such that $\{L_t^+\}_{0\leq t\leq T}\in{\cal{S}^{\mathrm{2}}}(0,T;{\mathbf{R}})$ and $L_T\leq\xi$. A solution
of this equation is a triple $(Y,Z,K)$ of progressively
measurable processes taking values in $\textbf{R}\times
{\textbf{R}}^d\times {\textbf{R}}_+$ and
satisfying
\begin{displaymath}
\begin{array}{rr}
\left\{\begin{array}{llll}
(i)\ \ Z\in{\cal{H}^{\mathrm{2}}}(0,T;{\mathbf{R}}^d),\ \ \ Y,\ K\in{\cal{S}^{\mathrm{2}}}(0,T;{\mathbf{R}})\\
(ii)\ \ Y_t=\xi +\int_t^Tg\left(s,Y_s,Z_s\right)
ds+K_T-K_t-\int_t^TZ_s\cdot dB_s,\ \ \forall t\in[0,T].\\
(iii)\ \ P-a.s,\ L_t\leq Y_t,\ \forall t\in[0,T]$ and
$\int_0^T(Y_t-L_t)dK_t=0.\\
(iv)\ \ K$ is continuous and increasing, $K_0=0.
\end{array}\right.
\end{array}
\end{displaymath}

The RBSDEs in Definiton 2.1 is called RBSDEs with parameter $(g,T,\xi,L)$, which is introduced in El Karoui et al. [8]. By Matoussi [11] or Hamad$\grave{\mathrm{e}}$ne et al. [12], we can get, under assumptions (A1) and (A2),  the RBSDEs has at least one solution. In particular, it has a minimal solution $(\underline{Y},\underline{Z},\underline{K})$ and a maximal solution $(\overline{{Y}},\overline{Z},\overline{K})$ in the sense that, for any solution  $({Y},{Z},{K})$ of this equation, we have $P-a.s,\ \underline{Y}\leq Y\leq\overline{{Y}}$ and $\underline{K}\geq K\geq\overline{{K}}.$ Note that in the following, we always assume that $\{L_t\}_{0\leq t\leq T},$ is a continuous progressively measurable real-valued process such that $\{L_t^+\}_{0\leq t\leq T}\in{\cal{S}^{\mathrm{2}}}(0,T;{\mathbf{R}}).$

\textbf{Remark 1} By Definition 2.1, we get that the solution of RBSDEs with parameter $(g,T,\xi,L)$ is only dependent on generator $g(t,y,z)|_{(t,y,z)\in[0,T]\times[L_t,+\infty)\times{\mathbf{R}}^{\mathit{d}}},$ not dependent on $g(t,y,z)|_{(t,y,z)\in[0,T]\times(-\infty, L_t[\times{\mathbf{R}}^{\mathit{d}}}.$ In other words, if $P-a.s., $
$$g_1(t,y,z)=g_2(t,y,z),\ \ \forall(t,y,z)\in[0,T]\times[L_t,+\infty)\times{\mathbf{R}}^{\mathit{d}},$$
then RBSDE with parameters $(g_1,T,\xi,L)$ and $(g_2,T,\xi,L)$ both have (or do not have) solutions. If one has solution(s), then another has the same solution(s).

The following Lemma 2.1 gives a priori estimation for RBSDEs, under assumptions (A1) and (A2), which can be proved using a standard argument given in Briand et al [2, Proposition 2.2] for BSDEs. We omit its proof.

\textbf{Lemma 2.1} Let $g$ satisfy the assumptions (A1) and (A2), $\xi\in L^2({\mathcal {F}}_T),$ $(Y,Z,K)$ be an arbitrary solution of RBSDEs with parameter $(g,\xi,T,L)$. Then there exists a constant $C$ depending only on $\lambda$ in (A1) and $T$,
such that for two stopping times $\sigma,\tau$ satisfying $0\leq \sigma<\tau\leq T,$ we have
\begin{eqnarray*}
&&E\left[\sup_{s\in [\sigma,\tau]}|{Y}_s|^2+\int_{\sigma}^\tau|{Z}_s|^2ds+|K_\tau-K_\sigma|^2|{\cal{F}}_\sigma\right]\\
&\leq& CE\left[|{Y}_\tau|^2+\left(\int_{\sigma}^\tau\gamma_sds\right)^2+\sup_{s\in[\sigma,\tau]}(L_s^+)^2|{\cal{F}}_\sigma\right].
\end{eqnarray*}

Now, we introduce a stochastic differential equation (SDE). Suppose $b(\cdot,\cdot,\cdot):\Omega\times [0,T]\times \textbf{R}^n\mapsto \textbf{R}^n$ and $\sigma(\cdot,\cdot,\cdot):\Omega\times [0,T]\times \textbf{R}^n\mapsto \textbf{R}^{n\times d}$ and always satisfy the following two conditions (H1) and (H2) in this paper.

(H1) (Lipschitz condition) there exists a constant $\mu>0$ such that $P$-$a.s.$, for each $t\in[0,T]$ and $x,y\in \textbf{R}^{\mathit{n}},$
$$|{b}(t,x)-b(t,y) |+|\sigma(t,x)-\sigma(t,y) |\leq \mu|x-y|.$$

(H2) (Linear-growth) there exists a constant $\nu>0$ such that for each $t\in[0,T]$ and $x\in \textbf{R}^{\mathit{n}},$
$$|b(t,x)|+|\sigma(t,x)|\leq \nu\left(1+|x|\right).$$

Given $(t,x)\in[0,T[\times \textbf{R}^n$, by SDE theory, the following SDE:
$$\left\{
  \begin{array}{ll}
     X_s^{t,x}=x +\int_t^sb(u,X_u^{t,x})du+\int_t^s\sigma(u,X_u^{t,x})dB_u,\ \ \  s\in]t,T],\\
     X_s^{t,x}=x\ \ \   s\in[0,t],
  \end{array}
\right.$$
has a unique continuous adapted solution $(X_s^{t,x})_{s\geq0}$.

\textbf{Remark 2} From the classical SDE theory, we know $E\left[|X_s^{t,x}-x|^p\right]$ is continuous in $s$ and $E\left[\sup_{0\leq s\leq T}|X_s^{t,x}|^p\right]<\infty,$ for each $p\geq 1.$

\textbf{Lemma 2.2 }  (Hewitt and Stromberg [13, Lemma 18.4]) Let $f$ be a Lebesgue integrable function on the interval $[0,T]$. Then for almost every $t\in[0,T[$, we have
$$\lim\limits_{\varepsilon\rightarrow0^+}\frac{1}{\varepsilon}
\int_{t}^{t+\varepsilon}|f(u)-f(t)|ds=0.$$

\textbf{Lemma 2.3 }  (Fan and Jiang [6, Corollary 1])  Let $g$ satisfy the assumptions (A1) and (A2), $ (\eta,x,q)\in{\cal{S}^{\mathrm{2}}}(0,T;{\mathbf{R}})\times {\mathbf{R}}^{\mathit{n}}\times {\mathbf{R}}^{\mathit{n}}$. Then there exists a non-negative process sequence $\{(\psi^m_t)_{t\in[0,T]}\}_{m=1}^{\infty}\in{\cal{H}^{\mathrm{2}}}(0,T;{\mathbf{R}})$ depending on $(\eta,x,q)$ such that $\lim\limits_{m\rightarrow\infty}\|\psi^m_t\|_{\cal{H}^{\mathrm{2}}}=0$, and $dP\times dt-a.s.,$ for any $m\in \textbf{N}$ and $(\bar{y},\bar{z},\bar{x})\in {\mathbf{R}}^{1+d+n}$,
$$|g(t,\bar{y},\bar{z}+\sigma^\ast(t,\bar{x})q)-g(t,\eta_t,\sigma^\ast(t,x)q)|\leq 2(m+\kappa)\left(|\bar{y}-\eta_t|+|\bar{z}|+|\bar{x}-x|\right)+\psi^m_t,$$
where $\kappa=\lambda(1+|q|\mu),$ $\lambda$ is the constant in (A1) and $\mu$ is the constant in (H1).

By the proof of Jiang [4, Proposition 2.2], we can get the following Lemma 2.4.

\textbf{Lemma 2.4 }  Let $q>p\geq1$. Let $(\psi_t)_{t\in[0,T]}$ be a real-valued,
progressively measurable process and
$E\left[\int_0^T|\psi_t|^qdt<\infty\right]$, then for almost every $t\in[0,T[$ and any stopping time $\tau\in]0,T-t],$ we have
$$\psi_t=\lim\limits_{\varepsilon\rightarrow0^+}\frac{1}{\varepsilon}
\int_{t}^{t+\varepsilon\wedge\tau}\psi_sds.$$
in $L^p({\mathcal {F}}_T)$ sense.
\section{Representation theorem for RBSDEs}
The following is a representation theorem for generators of RBSDEs, which is the main result of this paper.

\textbf{Theorem 3.1 }  Let $g$ satisfy the assumptions (A1) and (A2). Then for each $\eta\in {\cal{S}^{\mathrm{2}}}(0,T;{\mathbf{R}})$ satisfying $\eta_t> L_t$ for each $t\in[0,T]$, each $(x,q)\in{\mathbf{R}}^{\mathit{n}}\times {\mathbf{R}}^{\mathit{n}}$ and almost every $t\in[0,T[$, there exists a stopping time $\tau \in]0,T-t]$ depending on $(t,\eta_t,x,q)$, such that
$$g\left(t,\eta_t,\sigma^\ast(t,x)q\right)+q\cdot b(t,x)=\lim\limits_{\varepsilon\rightarrow0^+}\frac{1}{\varepsilon}
\left(Y_t^{t+\varepsilon\wedge\tau}-\eta_t-E\left[K_{t+\varepsilon\wedge\tau}^{t+\varepsilon\wedge\tau}
-K_t^{t+\varepsilon\wedge\tau}|{\cal{F}}_t\right]\right),\eqno(1)$$
in $L^p({\mathcal {F}}_T)$ sense for $1\leq p<2,$ where $(Y_s^{t+\varepsilon\wedge\tau},Z_s^{t+\varepsilon\wedge\tau},K_s^{t+\varepsilon\wedge\tau})$ is an arbitrary solution of RBSDEs with parameter $(g,t+\varepsilon\wedge\tau,\eta_{t}+q\cdot(X_{t+\varepsilon\wedge\tau}^{t,x}-x),L).$

\textit{Proof.} For $\eta\in {\cal{S}^{\mathrm{2}}}(0,T;{\mathbf{R}})$ satisfying $\eta_t> L_t$ for each $t\in[0,T]$, and $(t,x,q)\in[0,T[\times{\mathbf{R}}^{\mathit{n}}\times {\mathbf{R}}^{\mathit{n}}$, we define the following stopping time:
$$\tau:=\inf\left\{s\geq0:\eta_{t}+q\cdot(X_{t+s}^{t,x}-x)\leq L_{t+s}\right\}\wedge (T-t).\eqno(2)$$
By the fact $\eta_t> L_t$ for each $t\in[0,T]$, and the continuity of $X_{t+s}^{t,x}$, we have
$0<\tau\leq T-t$ and $$\eta_{t}+q\cdot(X_{t+{s\wedge\tau}}^{t,x}-x)\geq L_{t+{s\wedge\tau}}, \ \ \forall s\in[0,T].\eqno(3)$$
For $\varepsilon\in]0,T-t]$, let $(Y_s^{t+\varepsilon\wedge\tau},Z_s^{t+\varepsilon\wedge\tau},K_s^{t+\varepsilon\wedge\tau})$ be a solution of RBSDE with parameter $(g,t+\varepsilon\wedge\tau,\eta_{t}+q\cdot(X_{t+\varepsilon\wedge\tau}^{t,x}-x),L)$ and set
\begin{eqnarray*}
\ \ \ \ \ \ \tilde{Y}_s^{{t+\varepsilon\wedge\tau}}&:=&Y_s^{{t+\varepsilon\wedge\tau}}-(\eta_t+q\cdot(X_{s}^{t,x}-x)),\ \ \ \ \ \ \ \ s\in[t,{t+\varepsilon\wedge\tau}], \ \ \ \ \  \ \ \ \ \ \ \ \ (4)\\
\tilde{Z}_s^{t+\varepsilon\wedge\tau}&:=&Z_s^{t+\varepsilon\wedge\tau}-\sigma^\ast(s,X_{s}^{t,x})q,\ \ \ \ \ \ \ \ \ \ \ \ \ s\in[t,{t+\varepsilon\wedge\tau}],\\
\tilde{L}_s&:=&L_{t+\varepsilon\wedge\tau}-(\eta_t+q\cdot(X_{t+\varepsilon\wedge\tau}^{t,x}-x)),\ \ \ \ \ \ \ \ s\in]t+\varepsilon\wedge\tau,T],\\
\tilde{L}_s&:=&L_s-(\eta_t+q\cdot(X_{s}^{t,x}-x)), \ \ \ \ \ \ \ \ \ \ \ s\in[t,{t+\varepsilon\wedge\tau}],\\
\tilde{L}_s&:=&L_s-(\eta_s+q\cdot(X_{s}^{t,x}-x))=L_s-\eta_s, \ \ \ \ \ \ \ \ \ \ \ \ \ s\in[0,t[.
\end{eqnarray*}
Then by the fact $\eta_t> L_t$ for each $t\in[0,T],$ and (3), we have
$$\tilde{L}_s\leq0, \ \ \ \ \ s\in[0,T].\eqno(5)$$
Applying It\^{o}'s formula to $\tilde{Y}_s^{t+\varepsilon\wedge\tau},$ for $s\in[t,{t+\varepsilon\wedge\tau}],$ we have
\begin{eqnarray*}\ \tilde{Y}_s^{t+\varepsilon\wedge\tau}&=&\int_{s}^{{t+\varepsilon\wedge\tau}}g(r,\tilde{Y}_r^{t+\varepsilon\wedge\tau}+\eta_t
+q\cdot(X_{r}^{t,x}-x),
\tilde{Z}_r^{t+\varepsilon\wedge\tau}+\sigma^\ast(r,X_{r}^{t,x})q)dr\\
 &&+\int_{s}^{{t+\varepsilon\wedge\tau}}q\cdot b(r,X_{r}^{t,x})dr+{K}_{t+\varepsilon\wedge\tau}^{t+\varepsilon\wedge\tau}-{K}_{s}^{t+\varepsilon\wedge\tau}
 -\int_{s}^{t+\varepsilon\wedge\tau}\tilde{Z}_r^{t+\varepsilon\wedge\tau}dB_r. \ \ \ (6)
\end{eqnarray*}
Set
$$\tilde{g}(r,\tilde{y},\tilde{z})
:=\left\{
 \begin{array}{ll}
  0, \ \ \ \ \ \ \  r\in]{t+\varepsilon\wedge\tau},T],\\
  g(r,\tilde{y}+\eta_t+q\cdot(X_{r}^{t,x}-x),\tilde{z}+\sigma^\ast(r,X_{r}^{t,x})q)\\ \ \ \ \ +q\cdot b(r,X_{r}^{t,x}), \ \ \ \ \ \ \ \ \ \ r\in[t,{t+\varepsilon\wedge\tau}],\\
  g(r,\tilde{y}+\eta_r,\tilde{z}+\sigma^\ast(r,x)q)+q\cdot b(r,x),\ \ \ \ \ \ \ \ \  r\in[0,t[.
  \end{array}
  \right.
$$
By (A1) and (H2), we have
\begin{eqnarray*}\
\ \ \ \tilde{g}(r,\tilde{y},\tilde{z})&\leq& \lambda(\gamma_r+|\tilde{y}|+\sup_{u\in[0,T]}|\eta_u|+|q\cdot(X_{r}^{t,x}-x)|+|\tilde{z}+\sigma^\ast(r,X_{r}^{t,x})q|)\\
&&+|q\cdot b(r,X_{r}^{t,x})|\\
&\leq&(1+\lambda)\left(\tilde{\gamma}_r+|\tilde{y}|+|\tilde{z}|\right)\ \ \ \ \ \ \ \ \ \ \ \ \ \ \ \ \ \ \ \ \ \ \ \ \ \ \ \ \ \ \ \ \ \ \ \ \ \ \ \ \ \ \ \ \ \ \ \ \ \ \ \ \ (7)
\end{eqnarray*}
where $\tilde{\gamma}_r=(\gamma_r+\sup_{u\in[0,T]}|\eta_u|+|q|(|X_{r}^{t,x}|+|x|))+2\nu|q|(1+|X_{r}^{t,x}|).$ By Remark 2, $\eta\in {\cal{S}^{\mathrm{2}}}(0,T;{\mathbf{R}})$ and  $\gamma_t\in{\cal{H}^{\mathrm{2}}}(0,T;{\mathbf{R}}),$ we get $\tilde{\gamma}_t\in{\cal{H}^{\mathrm{2}}}(0,T;{\mathbf{R}})$. Then by (7) and (A2), we can check that $\tilde{g}$ also satisfy the assumptions (A1) and (A2). Let $({\hat{y}}_s^{t},{\hat{z}}_s^{t},{\hat{k}}_s^{t})$ be a solution of RBSDEs with parameter $(\tilde{g},t,\tilde{Y}_t^{t+\varepsilon\wedge\tau},\tilde{L}).$ By (6), it is not difficult to check that there exists a solution $({y}_s^{t+\varepsilon\wedge\tau},{z}_s^{t+\varepsilon\wedge\tau},{k}_s^{t+\varepsilon\wedge\tau})$ of RBSDEs with parameter $(\tilde{g},T,0,\tilde{L})$ such that, for $s\in[t,{t+\varepsilon\wedge\tau}],$
$$y_s^{{t+\varepsilon\wedge\tau}}=\tilde{Y}_s^{{t+\varepsilon\wedge\tau}}, \ z_s^{t+\varepsilon\wedge\tau}=\tilde{Z}_s^{t+\varepsilon\wedge\tau}, \ k_s^{t+\varepsilon\wedge\tau}={K}_s^{t+\varepsilon\wedge\tau}-{K}_t^{t+\varepsilon\wedge\tau}+\hat{k}_t^t,\eqno(8)$$
for $s\in]{t+\varepsilon\wedge\tau},T],$
$$y_s^{{t+\varepsilon\wedge\tau}}=0, \ z_s^{t+\varepsilon\wedge\tau}=0, \ k_s^{t+\varepsilon\wedge\tau}={K}_{t+\varepsilon\wedge\tau}^{t+\varepsilon\wedge\tau}
-{K}_t^{t+\varepsilon\wedge\tau}+\hat{k}_t^t,\eqno(9)$$
and for $s\in[0,t[,$
$$y_s^{{t+\varepsilon\wedge\tau}}=\hat{y}_s^{t}, \  z_s^{t+\varepsilon\wedge\tau}=\hat{z}_s^t,\ k_s^{t+\varepsilon\wedge\tau}=\hat{k}_s^t.\eqno(10)$$
Therefore by (5)-(10) and Lemma 2.1, we get there exists a constant $\tilde{C}$ depending only on $T$ and $\lambda,$
such that
$$E\left[\sup_{s\in[t,{t+\varepsilon\wedge\tau}]}|\tilde{Y}_s^{t+\varepsilon\wedge\tau}|^2
+\int_{t}^{t+\varepsilon\wedge\tau}|\tilde{Z}_s^{t+\varepsilon\wedge\tau}|^2ds\right]\leq \tilde{C}E\left[\left(\int_{t}^{t+\varepsilon\wedge\tau}\tilde{\gamma}_sds\right)^2\right].\ \ \ \ \eqno(11)$$
Set
\begin{eqnarray*}\
M^{\varepsilon,\tau}_t&:=&\frac{1}{\varepsilon}E\left[\int_{t}^{t+\varepsilon\wedge\tau}g(r,\tilde{Y}_r^{t+\varepsilon\wedge\tau}
+\eta_t+q\cdot(X_{r}^{t,x}-x),
\tilde{Z}_r^{t+\varepsilon\wedge\tau}+\sigma^\ast(r,X_{r}^{t,x})q)dr|{\cal{F}}_t\right]\\
P^{\varepsilon,\tau}_t&:=&\frac{1}{\varepsilon}E\left[\int_{t}^{t+\varepsilon\wedge\tau}g(r,\eta_r,
\sigma^\ast(r,x)q)dr|{\cal{F}}_t\right],\\
U^{\varepsilon,\tau}_t&:=&\frac{1}{\varepsilon}E\left[\int_{t}^{t+\varepsilon\wedge\tau}q\cdot b(r,X_{r}^{t,x})dr|{\cal{F}}_t\right],
\end{eqnarray*}
By (4), we have $Y_t^{t+\varepsilon\wedge\tau}-\eta_t=\tilde{Y}_t^{t+\varepsilon\wedge\tau}$. Then by (6), we have
\begin{eqnarray*} &&\frac{1}{\varepsilon}\left(Y_t^{t+\varepsilon\wedge\tau}-\eta_t-E\left[K_{t+\varepsilon\wedge\tau}^{t+\varepsilon\wedge\tau}-K_t^{t+\varepsilon\wedge\tau}|{\cal{F}}_t\right]\right)
-g(t,\eta_t,\sigma^\ast(t,x)q)-q\cdot b(t,x)\\
&=&\frac{1}{\varepsilon}\left(\tilde{Y}_t^{t+\varepsilon\wedge\tau}-E\left[{K}_{t+\varepsilon\wedge\tau}^{t+\varepsilon\wedge\tau}
-{K}_t^{t+\varepsilon\wedge\tau}|{\cal{F}}_t\right]\right)
-g(t,\eta_t,\sigma^\ast(t,x)q)-q\cdot b(t,x)\\
&=&M^{\varepsilon,\tau}_t+U^{\varepsilon,\tau}_t-g(t,\eta_t,\sigma^\ast(t,x)q)-q\cdot b(t,x)\\
&=&\left(M^{\varepsilon,\tau}_t-P^{\varepsilon,\tau}_t\right)
+\left(P^{\varepsilon,\tau}_t-g(t,\eta_t,\sigma^\ast(t,x)q)\right)
+\left(U^{\varepsilon,\tau}_t-q\cdot b(t,x)\right).\ \ \ \ \ \ \ \ \ \ \ \ \ \ \ (12)
\end{eqnarray*}
Thus, we only need prove that (12) converges to $0$ in $L^p({\mathcal {F}}_T)$ sense for $1\leq p<2.$

By Jensen's inequality, H\"{o}lder's inequality and Lemma 2.3, we get that there exists a non-negative process sequence $\{(\psi^m_t)_{t\in[0,T]}\}_{m=1}^{\infty}\in{\cal{H}^{\mathrm{2}}}(0,T;{\mathbf{R}})$ depending on $(\eta_t,x,q)$ and $\lim\limits_{m\rightarrow\infty}\|\psi^m_t\|_{\cal{H}^{\mathrm{2}}}=0,$ such that, for any $m\geq1,$
\begin{eqnarray*}
&&E|M^{\varepsilon,\tau}_t-P^{\varepsilon,\tau}_t|^2\\
&\leq&\frac{1}{\varepsilon^2}E\Bigg[\int_{t}^{t+\varepsilon\wedge\tau}|g(r,\tilde{Y}_r^{t+\varepsilon\wedge\tau}+\eta_t
+q\cdot(X_{r}^{t,x}-x),
\tilde{Z}_r^{t+\varepsilon\wedge\tau}+\sigma^\ast(r,X_{r}^{t,x})q)\\&&-g(r,\eta_r,
\sigma^\ast(r,x)q)|dr\Bigg]^2\\
&\leq&\frac{1}{\varepsilon}E\Bigg[\int_{t}^{t+\varepsilon\wedge\tau}|g(r,\tilde{Y}_r^{t+\varepsilon\wedge\tau}+\eta_t
+q\cdot(X_{r}^{t,x}-x),
\tilde{Z}_r^{t+\varepsilon\wedge\tau}+\sigma^\ast(r,X_{r}^{t,x})q)\\&&-g(r,\eta_r,
\sigma^\ast(r,x)q)|^2dr\Bigg]\\
&\leq&\frac{1}{\varepsilon}E\Bigg[\int_{t}^{t+\varepsilon\wedge\tau}(2(m+\kappa)(|\tilde{Y}_r^{t+\varepsilon\wedge\tau}|
+|\tilde{Z}_r^{t+\varepsilon\wedge\tau}|
+|\eta_r-\eta_t|\\&&+(|q|+1)|X_{r}^{t,x}-x|)+\psi_r^m)^2dr\Bigg]\\
&\leq&\frac{1}{\varepsilon}E\Bigg[\int_{t}^{t+\varepsilon\wedge\tau} 8(m+\kappa)^2\Bigg((|\tilde{Y}_r^{t+\varepsilon\wedge\tau}|
+|\tilde{Z}_r^{t+\varepsilon\wedge\tau}|
+|\eta_r-\eta_t|\\&&+(|q|+1)|X_{r}^{t,x}-x|\Bigg)^2+2|\psi_r^m|^2dr\Bigg]
\end{eqnarray*}
where $\kappa=\lambda(1+|q|\mu).$ Then by the above inequality, (11) and H\"{o}lder's inequality, we have for $m\geq1,$
\begin{eqnarray*}
&&E|M^{\varepsilon,\tau}_t-P^{\varepsilon,\tau}_t|^2\\
&\leq&32(m+\kappa)^2\tilde{C}\frac{1}{\varepsilon}E\left[\left(\int_{t}^{t+\varepsilon}\tilde{\gamma}_rdr\right)^2\right]
+32(m+\kappa)^2\frac{1}{\varepsilon}E\left[\int_{t}^{t+\varepsilon}|\eta_r-\eta_t|^2dr\right]\\
&&+32(m+\kappa)^2(|q|+1)^2\frac{1}{\varepsilon}E\left[\int_{t}^{t+\varepsilon}|X_{r}^{t,x}-x|^2dr\right]
+\frac{2}{\varepsilon}E\left[\int_{t}^{t+\varepsilon}|\psi_r^m|^2dr\right]\\
&\leq&32(m+\kappa)^2\tilde{C}E\left[\int_{t}^{t+\varepsilon}|\tilde{\gamma}_r|^2 dr\right]+32(m+\kappa)^2\frac{1}{\varepsilon}E\left[\int_{t}^{t+\varepsilon}|\eta_r-\eta_t|^2dr\right]\\
&&+32(m+\kappa)^2(|q|+1)^2\frac{1}{\varepsilon}E\left[\int_{t}^{t+\varepsilon}|X_{r}^{t,x}-x|^2dr\right]
+\frac{2}{\varepsilon}E\left[\int_{t}^{t+\varepsilon}|\psi_r^m|^2dr\right], \ \ \ \ (13)
\end{eqnarray*}
Since $\tilde{\gamma}_t\in{\cal{H}^{\mathrm{2}}}(0,T;{\mathbf{R}})$, then by Fubini's theorem and absolute continuity of integral, we have
$$\lim\limits_{\varepsilon\rightarrow0^+}E\left[\int_{t}^{t+\varepsilon}|\tilde{\gamma}_r|^2dr\right]
=\lim\limits_{\varepsilon\rightarrow0^+}\int_{t}^{t+\varepsilon}E|\tilde{\gamma}_r|^2dr=0.\eqno(14)$$
Since $\eta\in {\cal{S}^{\mathrm{2}}}(0,T;{\mathbf{R}})$, we can deduce $E|\eta_r-\eta_t|^2$ is continuous in $r$. Then by Fubini's Theorem, we have
$$\lim\limits_{\varepsilon\rightarrow0^+}\frac{1}{\varepsilon}E\left[\int_{t}^{t+\varepsilon}|\eta_r-\eta_t|^2dr\right]
=\lim\limits_{\varepsilon\rightarrow0^+}\frac{1}{\varepsilon}\int_{t}^{t+\varepsilon}E|\eta_r-\eta_t|^2dr=0.\eqno(15)$$
By Fubini's theorem and Remark 2, we have,
$$\lim\limits_{\varepsilon\rightarrow0^+}\frac{1}{\varepsilon}E\left[\int_{t}^{t+\varepsilon}|X_{r}^{t,x}-x|^2dr\right]
=\lim\limits_{\varepsilon\rightarrow0^+}\frac{1}{\varepsilon}\int_{t}^{t+\varepsilon}E|X_{r}^{t,x}-x|^2dr=0. \eqno(16)$$
By Fatou's lemma, Fubini's theorem and $\lim\limits_{m\rightarrow\infty}\|\psi^m_t\|_{\cal{H}^{\mathrm{2}}}=0,$ we have
\begin{eqnarray*}
\int_{0}^T\liminf_{m\rightarrow\infty}E|\psi_r^m|^2dr\leq \liminf_{m\rightarrow\infty}\int_{0}^TE|\psi_r^m|^2dr&=&
\liminf_{m\rightarrow\infty}E\int_{0}^T|\psi_r^m|^2dr\\
&=&\lim\limits_{n\rightarrow\infty}
\|\psi^m_t\|_{\cal{H}^{\mathrm{2}}}^2\\&=&0.
\end{eqnarray*}
Thus, for almost every $t\in[0,T[$, we have
$$\liminf_{m\rightarrow\infty}E|\psi_t^m|^2=0.\eqno(17)$$
Then by (13)-(16), Fubini's theorem, Lemma 2.2 and (17), we get that for almost every $t\in[0,T[$,
\begin{eqnarray*}\ \ \ \ \ \ \ \lim\limits_{\varepsilon\rightarrow0^+}E|M^{\varepsilon,\tau}_t
-P^{\varepsilon,\tau}_t|^2&\leq&\liminf_{m\rightarrow\infty}\lim\limits_{\varepsilon\rightarrow0^+}\frac{2}{\varepsilon}
E\left[\int_{t}^{t+\varepsilon}|\psi_r^m|^2dr\right]\\
&=&\liminf_{m\rightarrow\infty}\lim\limits_{\varepsilon\rightarrow0^+}\frac{2}{\varepsilon}
\int_{t}^{t+\varepsilon}E|\psi_r^m|^2dr\\
&=&2\liminf_{m\rightarrow\infty}E|\psi_t^m|^2\\
&=&0.\ \ \ \ \ \ \ \ \ \ \ \ \ \ \ \ \ \ \ \ \ \ \ \ \ \ \ \ \ \ \ \ \ \ \ \ \ \ \ \ \ \ \ \ \ \ \ \ \ \ \ \ \ \ \ \ \ \ \ \ (18)
\end{eqnarray*}
For $1\leq p<2,$ by Jensen's inequality and Lemma 2.4, we have, for almost every $t\in[0,T[$,
\begin{eqnarray*}
&&\lim\limits_{\varepsilon\rightarrow0^+}E|P^{\varepsilon,\tau}_t-g(t,\eta_t,\sigma^\ast(t,x)q)|^p\\
&=&\lim\limits_{\varepsilon\rightarrow0^+}E\left|E\left[\frac{1}{\varepsilon}\int_{t}^{t+\varepsilon\wedge\tau} (g(r,\eta_r,
\sigma^\ast(r,x)q)dr-g(t,\eta_t,\sigma^\ast(t,x)q)|{\cal{F}}_t\right]\right|^p\\
&\leq&\lim\limits_{\varepsilon\rightarrow0^+}E\left|\frac{1}{\varepsilon}\int_{t}^{t+\varepsilon\wedge\tau}g(r,\eta_r,
\sigma^\ast(r,x)q)dr-g(t,\eta_t,\sigma^\ast(t,x)q)\right|^p\\
&=&0.\ \ \ \ \ \ \ \ \ \ \ \ \ \ \ \ \ \ \ \ \ \ \ \  \ \ \ \ \ \ \ \ \ \ \ \ \ \ \ \ \ \ \ \ \ \ \ \ \ \ \ \ \ \ \ \ \ \ \ \ \ \ \ \ \ \ \ \ \ \ \ \ \ \ \ \ \ \ \ \ \ \ \ \ \ \ \ \ \ \ \ \ \ \ (19)
\end{eqnarray*}
For $1\leq p<2,$ by Jensen's inequality and Lemma 2.4,  for almost every $t\in[0,T[$,
we also have,
\begin{eqnarray*}
&&\lim\limits_{\varepsilon\rightarrow0^+}E|U^{\varepsilon,\tau}_t-q\cdot b(t,x)|^p\\
&=&\lim\limits_{\varepsilon\rightarrow0^+}E\left|E\left[\frac{1}{\varepsilon}\int_{t}^{t+\varepsilon\wedge\tau}
(q\cdot b(r,X_{r}^{t,x})-q\cdot b(r,x)+q\cdot b(r,x))dr-q\cdot b(t,x)|{\cal{F}}_t\right]\right|^p\\
&\leq&\lim\limits_{\varepsilon\rightarrow0^+}2^{p-1}E\left|E\left[\frac{1}{\varepsilon}\int_{t}^{t+\varepsilon\wedge\tau}(q\cdot b(r,X_{r}^{t,x})-q\cdot b(r,x))dr|{\cal{F}}_t\right]\right|^p\\
&&+\lim\limits_{\varepsilon\rightarrow0^+}2^{p-1}E\left|E\left[\frac{1}{\varepsilon}\int_{t}^{t+\varepsilon\wedge\tau}q\cdot b(r,x)dr-q\cdot b(t,x)|{\cal{F}}_t\right]\right|^p\\
&\leq&\lim\limits_{\varepsilon\rightarrow0^+}\frac{2^{p-1}}{\varepsilon^p}E\left(\int_{t}^{t+\varepsilon\wedge\tau}|q\cdot b(r,X_{r}^{t,x})-q\cdot b(r,x)|dr\right)^p\\
&&+\lim\limits_{\varepsilon\rightarrow0^+}2^{p-1}E\left|\frac{1}{\varepsilon}\int_{t}^{t+\varepsilon\wedge\tau}q\cdot b(r,x)dr-q\cdot b(t,x)\right|^p\\
&=&\lim\limits_{\varepsilon\rightarrow0^+}\frac{2^{p-1}}{\varepsilon^p}E\left(\int_{t}^{t+\varepsilon}|q\cdot b(r,X_{r}^{t,x})-q\cdot b(r,x)|dr\right)^p.
\end{eqnarray*}
By the above inequality, H\"{o}lder's inequality, Fubini's theorem, (H1) and Remark 2, we have, for almost every $t\in[0,T[$,
$$\lim\limits_{\varepsilon\rightarrow0^+}E|U^{\varepsilon,\tau}_t-q\cdot b(t,x)|^p
\leq\lim\limits_{\varepsilon\rightarrow0^+}\frac{2^{p-1}}{\varepsilon}|q|^p\mu^p\int_{t}^{t+\varepsilon}E|X_{r}^{t,x}-x|^pdr
=0.\eqno(20)$$
By (12) and (18)-(20), we obtain (1). The proof is complete. \ \ $\Box$

Now, we consider two special cases of Theorem 3.1. When the obstacle $L_t$ is a special regular process, we will have the following Corollary 3.2.

\textbf{Corollary 3.2 } Let $g$ satisfy assumptions (A1) and (A2), and $L_t$ is an It\^{o} process:
$$L_t=L_0+\int_0^tU_sds+\int_0^tV_sdB_s,\ \ t\in[0,T],$$
where $U_s$ and $V_s$ are two progressively measurable processes such that for each $t\in[0,T], P-a.s.,\ \int_0^t(|U_s|+|V_s|^2)ds<\infty$ and $g(t,L_t,V_t)+U_t\geq0.$ Then for each $\eta\in {\cal{S}^{\mathrm{2}}}(0,T;{\mathbf{R}})$ satisfying $\eta_t> L_t$ for each $t\in[0,T]$, each $  (x,q)\in{\mathbf{R}}^{\mathit{n}}\times {\mathbf{R}}^{\mathit{n}}$ and almost every $t\in[0,T[$, there exists a stopping time $\tau \in]0,T-t]$ depending on $(t,\eta_t,x,q)$, such that
$$g\left(t,\eta_t,\sigma^\ast(t,x)q\right)+q\cdot b(t,x)=\lim\limits_{\varepsilon\rightarrow0^+}\frac{1}{\varepsilon}
\left(Y_t^{t+\varepsilon\wedge\tau}-\eta_t\right),$$
in $L^p({\mathcal {F}}_T)$ sense for $1\leq p<2,$ where $(Y_s^{t+\varepsilon\wedge\tau},Z_s^{t+\varepsilon\wedge\tau},K_s^{t+\varepsilon\wedge\tau})$ is an arbitrary solution of RBSDEs with parameter $(g,t+\varepsilon\wedge\tau,\eta_{t}+q\cdot(X_{t+\varepsilon\wedge\tau}^{t,x}-x),L).$

\textit{Proof.}  Let $\tau$ be the stopping time defined in (2), $(Y_s^{t+\varepsilon\wedge\tau},Z_s^{t+\varepsilon\wedge\tau},K_s^{t+\varepsilon\wedge\tau})$ is an arbitrary solution of RBSDEs with parameter $(g,t+\varepsilon\wedge\tau,\eta_{t}+q\cdot(X_{t+\varepsilon\wedge\tau}^{t,x}-x),L).$ By El Karoui et al. [8, Proposition 4.2], we have
$$P-a.s.,\ \ 0\leq K_{t+\varepsilon\wedge\tau}^{t+\varepsilon\wedge\tau}-K_t^{t+\varepsilon\wedge\tau}
\leq\int_t^{t+\varepsilon\wedge\tau}1_{\{Y_s^{t+\varepsilon\wedge\tau}=L_s\}}(g(s,L_s,V_s)+U_s)^-ds=0.$$
From Theorem 3.1, the proof is complete. \ \ $\Box$

When the obstacle $L_t$ has a upper bound, we will have the following Corollary 3.3.

\textbf{Corollary 3.3 } Let $g$ satisfy assumptions (A1)-(A3), and there exists a constant $C$ such that $\sup_{t\in [0,T]}L_t\leq C.$ Then for each $(y,x,q)\in]C,+\infty)\times{\mathbf{R}}^{\mathit{n}}\times{\mathbf{R}}^{\mathit{n}}$ and for almost every $t\in[0,T[$, there exists a stopping time $\tau>0$ depending on $(t,y,x,q)$, such that
$$g(t,y,\sigma^\ast(t,x)q)+q\cdot b(t,x)=\lim\limits_{\varepsilon\rightarrow0^+}\frac{1}{\varepsilon}
\left(Y_t^{t+\varepsilon\wedge\tau}-y\right),\eqno(21)$$
in $L^p({\mathcal {F}}_T)$ sense for $1\leq p<2,$ where $(Y_s^{t+\varepsilon\wedge\tau},Z_s^{t+\varepsilon\wedge\tau},K_s^{t+\varepsilon\wedge\tau})$ is an arbitrary solution of RBSDEs with parameter $(g,t+\varepsilon\wedge\tau,y+q\cdot(X_{t+\varepsilon\wedge\tau}^{t,x}-x),L)$ satisfying $P-a.s.,\ Y_s^{t+\varepsilon\wedge\tau}\geq C$ for each $s\in[0,t+\varepsilon\wedge\tau].$

\textit{Proof.} For each $(t,x,q)\in[0,T]\times{\mathbf{R}}^{\mathit{d}}\times {\mathbf{R}}^{\mathit{n}}$ and $y> C,$ we define the following stopping time:
 $$\tau=\inf\left\{s\geq0:y+q\cdot(X_{t+s}^{t,x}-x)\leq C\right\}\wedge (T-t),$$
By $y> C$ and the continuity of $X_{t+s}^{t,x}$, we have
$0<\tau\leq T-t$ and $$y+q\cdot(X_{t+{s\wedge\tau}}^{t,x}-x)\geq C\geq L_{t+{s\wedge\tau}}, \ \ \forall s\in[0,T].\eqno(22)$$
For $\varepsilon\in]0,T-t]$, let $(y_s^{C,t+\varepsilon\wedge\tau},z_s^{C,t+\varepsilon\wedge\tau},k_s^{C,t+\varepsilon\wedge\tau})$ be a solution of RBSDEs with parameter $(g,t+\varepsilon\wedge\tau,y+q\cdot(X_{t+{s\wedge\tau}}^{t,x}-x),C)$. Since the obstacle is $C,$ by El Karoui et al. [8, Proposition 4.2] and (A3), we can get for each $s\in[0,t+\varepsilon\wedge\tau],$ $$ P-a.s.,\ k_s^{C,t+\varepsilon\wedge\tau}=0.\eqno(23)$$
Let $(Y_s^{t+\varepsilon\wedge\tau},Z_s^{t+\varepsilon\wedge\tau},K_s^{t+\varepsilon\wedge\tau})$ is a solution of RBSDEs with parameter $(g,t+\varepsilon\wedge\tau,y+q\cdot(X_{t+\varepsilon\wedge\tau}^{t,x}-x),L),$ such that for each $s\in[0,t+\varepsilon\wedge\tau], P-a.s.,\ Y_s^{t+\varepsilon\wedge\tau}\geq C.$ Clearly, we can check $(Y_s^{t+\varepsilon\wedge\tau},Z_s^{t+\varepsilon\wedge\tau},K_s^{t+\varepsilon\wedge\tau})$ is also a solution of RBSDEs with parameter $(g,t+\varepsilon\wedge\tau,y+q\cdot(X_{t+{s\wedge\tau}}^{t,x}-x),C)$. By (23) and Theorem 3.1, we can get (21).

We will show there exists a solution of RBSDEs with parameter $(g,t+\varepsilon\wedge\tau,y+q\cdot(X_{t+\varepsilon\wedge\tau}^{t,x}-x),L),$ such that for each $s\in[0,t+\varepsilon\wedge\tau], P-a.s.,\ Y_s^{t+\varepsilon\wedge\tau}\geq C.$  For $\varepsilon\in]0,T-t]$, we denote the maximal solution of RBSDEs with parameter $(g,t+\varepsilon\wedge\tau,y+q\cdot(X_{t+\varepsilon\wedge\tau}^{t,x}-x),L)$ and $(g,t+\varepsilon\wedge\tau,C,L)$ by $(\overline{Y}_s^{t+\varepsilon\wedge\tau},\overline{Z}_s^{t+\varepsilon\wedge\tau},\overline{K}_s^{t+\varepsilon\wedge\tau})$ and $(\overline{Y}_s^{C,t+\varepsilon\wedge\tau},\overline{Z}_s^{C,t+\varepsilon\wedge\tau},\overline{K}_s^{C,t+\varepsilon\wedge\tau})$, respectively. By (A3) and $\sup_{t\in [0,T]}L_t\leq C$, we can check that $(\hat{Y}_s^{C,t+\varepsilon\wedge\tau},\hat{Z}_s^{C,t+\varepsilon\wedge\tau},\hat{K}_s^{C,t+\varepsilon\wedge\tau})=(C,0,0),\ s\in[0,t+\varepsilon\wedge\tau],$ is a solution of RBSDEs with parameter $(g,t+\varepsilon\wedge\tau,C,L)$. By (22) and comparison theorem of RBSDEs with continuous and linear-growth generators (which can be proved by Hamad$\grave{\mathrm{e}}$ne et al. [12, Proposition 2.3] and a simple argument as the corresponding proof for BSDEs case), we have
$$\overline{Y}_s^{t+\varepsilon\wedge\tau}\geq \overline{Y}_s^{C,t+\varepsilon\wedge\tau}\geq \hat{Y}_s^{C,t+\varepsilon\wedge\tau}=C,\ \ s\in[0,t+\varepsilon\wedge\tau].$$
The proof is complete.\ \  $\Box$

Let $n=d,\ q=z,\ b(t,x)=0,\ \sigma(t,x)=1,\ x=0$ in Theorem 3.1. Then we have the following Corollary 3.4, immediately.
 
\textbf{Corollary 3.4 } Let $g$ satisfy the assumptions (A1) and (A2). Then for each $\eta\in {\cal{S}^{\mathrm{2}}}(0,T;{\mathbf{R}})$ satisfying $\eta_t> L_t$ for each $t\in[0,T],$ each $  z\in{\mathbf{R}}^{\mathit{d}}$ and almost every $t\in[0,T[$, there exists a stopping time $\tau>0$ depending on $(t,\eta_t,z)$, such that
$$g(t,\eta_t,z)=\lim\limits_{\varepsilon\rightarrow0^+}\frac{1}{\varepsilon}
\left(Y_t^{t+\varepsilon\wedge\tau}-\eta_t-E\left[K_{t+\varepsilon\wedge\tau}^{t+\varepsilon\wedge\tau}
-K_t^{t+\varepsilon\wedge\tau}|{\cal{F}}_t\right]\right),$$
in $L^p({\mathcal {F}}_T)$ sense for $1\leq p<2,$ where $(Y_s^{t+\varepsilon\wedge\tau},Z_s^{t+\varepsilon\wedge\tau},K_s^{t+\varepsilon\wedge\tau})$ is an arbitrary solution of RBSDE with parameter $(g,t+\varepsilon\wedge\tau,\eta_{t}+z\cdot(B_{t+\varepsilon\wedge\tau}-B_t),L).$

\textbf{Remark 3}
\begin{enumerate}
  \item From the proof of Theorem 3.1 and Corollary 3.2, we can get that the stopping time $\tau$ in the above representation theorems can be replaced by any stopping time $\sigma\in]0,\tau].$
  \item Our representation theorems for RBSDEs are obtained in local space. In fact, from Remark 1, it follows that it may be not true in whole space. Due to this fact, we only can use our representation theorem to study the local properties of generator of RBSDE. This is different from BSDEs case.

  \item The representation theorem for BSDEs is established in Jiang [3] under the additional condition that ${b}(t,x)$ and $\sigma(t,x)$ in SDE are both right continuous in $t$ and established in Jia [5] under the additional condition that ${b}(t,x)$ and $\sigma(t,x)$ in SDE are both independent on $t$. One can see that such conditions are eliminated in our representation theorems due to the uses of Lemma 2.2 and Lemma 2.4.
  \item Let $(Y,Z,K)$ be a solution of RBSDEs with parameter $(g,T,\xi,L)$. Indeed, if $L\equiv-\infty,$ then $K\equiv0.$ In this case, we can see RBSDEs will become BSDEs, and Theorem 3.1 will become a representation theorem for BSDEs with continuous and linear growth generators, which has been studied in Jia [5] and Fan and Jiang [6].
\end{enumerate}
\section{Some applications}
With the help of the representation theorem of RBSDEs, we can establish a general converse comparison theorems for RBSDEs.

\textbf{Theorem 4.1 }  Let generators $g_1$ and $g_2$ satisfy assumptions (A1) and (A2). If for any stopping time $\tau \in]0,T]$ and $\xi\in L^2({\mathcal{F}}_\tau)$ satisfying $\xi\geq L_\tau$, there exists a solution $(Y_t^{\tau,i},Z_t^{\tau,i},K_t^{\tau,i})$ of RBSDEs with parameter $(g_i,\tau,\xi,L),\ i=1,2,$  such that for each $t\in[0,T],$
$$P-a.s.,\  Y_{t\wedge\tau}^{\tau,1}\geq Y_{t\wedge\tau}^{\tau,2},\ \ \textmd{and}\ \ E[(K_\tau^{\tau,1}-K_{t\wedge\tau}^{\tau,1})|{\cal{F}}_t]\leq E[(K_\tau^{\tau,2}-K_{t\wedge\tau}^{\tau,2})|{\cal{F}}_t], \eqno(24)$$
then for each $\eta\in {\cal{S}^{\mathrm{2}}}(0,T;{\mathbf{R}})$ satisfying $\eta\geq L$, each $ z\in{\mathbf{R}}^{\mathit{d}}$, and almost every $t\in[0,T[$, we have
$$P-a.s.,\  \ g_1(t,\eta_t,z)\geq g_2(t,\eta_t,z).$$

\textit{Proof.} By Corollary 3.4, for each $\eta\in {\cal{S}^{\mathrm{2}}}(0,T;{\mathbf{R}})$ satisfying $\eta_t> L_t$ for each $t\in[0,T]$, each $z\in{\mathbf{R}}^{\mathit{d}}$ and almost every $t\in[0,T[$, there
exists a stopping time $\delta>0$ depending on $(t,\eta_t,z)$ and a subsequence $\{n_k\}_{k\geq1}$ of $\{n\}$, such that $P-a.s.,$
$$g_i(t,\eta_t,z)=\lim\limits_{k\rightarrow+\infty}n_k
\left(Y_t^{t+n_k^{-1}\wedge\delta,i}-\eta_t-E\left[\left(K_{t+n_k^{-1}\wedge\delta}^{t+n_k^{-1}\wedge\delta,i}
-K_t^{t+n_k^{-1}\wedge\delta,i}\right)|{\cal{F}}_t\right]\right), \eqno(25)$$
where $(Y_t^{t+n_k^{-1}\wedge\delta,i},Z_t^{t+n_k^{-1}\wedge\delta,i},K_t^{t+n_k^{-1}\wedge\delta,i})$ is an arbitrary solution of RBSDE with parameter $(g_i,t+n_k^{-1}\wedge\delta,\eta_{t}+z\cdot(B_{t+n_k^{-1}\wedge\delta}-B_t),L),\ i=1,2.$
By (24),  (25) and (A2), the proof is complete.\ \ $\Box$

Similarly, we can obtain the following Corrollary 4.2 from Corollary 3.3.

\textbf{Corollary 4.2 }  Let generators $g_1$ and $g_2$ satisfy assumptions (A1)-(A3), and there exists a constant $C$ such that $\sup_{t\in [0,T]}L_t\leq C$. If for any stopping time $\tau \in]0,T]$ and $\xi\in L^2({\mathcal{F}}_\tau),$ there exists a solution $(Y_t^{\tau,i},Z_t^{\tau,i},K_t^{\tau,i})$ of RBSDEs with parameter $(g_i,\tau,\xi,L)$ such that for each $t\in[0,T], P-a.s.,\ Y_{t\wedge\tau}^{\tau,i}\geq C,\ i=1,2,$ and
$$P-a.s.,\  Y_{t\wedge\tau}^{\tau,1}\geq Y_{t\wedge\tau}^{\tau,2},$$
then for each $(y,z)\in[C,+\infty)\times{\mathbf{R}}^{\mathit{d}}$ and almost every $t\in[0,T[$, we have
$$P-a.s.,\  \ g_1(t,y,z)\geq g_2(t,y,z).$$

\textbf{Remark 4}
\begin{enumerate}
\item Theorem 4.1 and Corollary 4.2 are both established for generators which are continuous with linear-growth in $(y,z)$, while converse comparison theorems of RBSDEs in Li and Tang [9] and Li and Gu [10] are  obtained under the conditions that generators satisfy Lipschitz condition and are continuous in $t$.
\item Theorem 3.4 and Corollary 3.5 both show that generators can be compared in local space. In fact, by Remark 1, we get that the generators can not be compared in whole space. This is different from BSDEs case.
\end{enumerate}

In the following, we will discuss some properties of RBSDEs, which have been considered in Jia [5], Fan and Jiang [6], and Fan et al. [7] in BSDEs case.

\textbf{Proposition 4.1} (Self-financing condition) Let $g$ satisfy assumptions (A1) and (A2),  and $\sup_{t\in [0,T]}L_t<0.$ Then the following two statements are equivalent:

(i) For almost every $t\in[0,T[$,
$$P-a.s.,\ \ g(t,0,0)=0;$$

(ii) There exists a solution $(Y_t,Z_t,K_t)$ of RBSDEs with parameter $(g,T,0,L)$ such that
$$P-a.s.,\ \ Y_t=0,\ \forall t\in[0,T].$$

\emph{Proof.} We can easily check that (i) implies (ii) by setting $(Y_t,Z_t,K_t)=(0,0,0)$. If (ii) holds, then we can check that for any stopping time $\tau\in]0,T]$, there exist a solution $(y_t^\tau,z_t^\tau,k_t^\tau)$ of RBSDEs with parameter $(g,\tau,0,L)$ such that for each $t\in[0,\tau],$
$$P-a.s.,\ \ y_t^\tau=0,\ \ k_t^\tau=0.$$
Then by Corollary 3.4, the proof is complete.\ \ \ $\Box$

\textbf{Proposition 4.2} (Zero-interest condition) Let $g$ satisfy assumptions (A1) and (A2), and there exists a constant $C$ such that $\sup_{t\in [0,T]}L_t\leq C.$ Then the following two statements are equivalent:

(i) For almost every $t\in[0,T[$ and each $y\in[C,+\infty),$
$$P-a.s.,\ \ g(t,y,0)=0;$$

(ii) For any $y\geq C$, there exists a solution $(Y_t,Z_t,K_t)$ of RBSDEs with parameter $(g,T,y,L)$ such that
$$P-a.s.,\ \ Y_t=y,\ \forall t\in[0,T].$$

\emph{Proof.} We can easily check that (i) implies (ii) by setting $(Y_t,Z_t,K_t)=(y,0,0)$. If (ii) holds, then we can check that for $y>C$, any stopping time $\tau\in]0,T]$, there exist a solution $(y_t^\tau,z_t^\tau,k_t^\tau)$ of RBSDEs with parameter $(g,\tau,y,L)$ such that for each $t\in[0,\tau],$
$$P-a.s.,\ \ y_t^\tau=y,\ \ k_t^\tau=0.$$
Then by Corollary 3.4 and (A2), the proof is complete.\ \ \ $\Box$

More generally, we have

\textbf{Proposition 4.3} Let $g$ satisfy assumptions (A1) and (A2). Then the following two statements are equivalent:

(i) For each $\eta\in {\cal{S}^{\mathrm{2}}}(0,T;{\mathbf{R}})$ satisfying $\eta\geq L,$ and almost every $t\in[0,T[$,
$$P-a.s.,\ \ g(t,\eta_t,0)=0;$$

(ii) For each $\eta\in {\cal{S}^{\mathrm{2}}}(0,T;{\mathbf{R}})$ satisfying $\eta\geq L,$ and each $t\in[0,T[$, if there exists a stopping time $\sigma_t\in]t,T]$ such that $\eta_t\geq L_s$ for each $s\in[t,\sigma_t],$ then there exists a solution $(Y_t,Z_t,K_t)$ of RBSDE with parameter $(g,\sigma_t,\eta_t,L)$ such that for each $s\in[t,\sigma_t],$
$$P-a.s.,\ \ Y_s=\eta_t.$$

\emph{Proof.} We firstly prove (i) implies (ii). Let $\eta\in {\cal{S}^{\mathrm{2}}}(0,T;{\mathbf{R}})$ satisfying $\eta\geq L,$ and for each $t\in[0,T[$ and there exists a stopping time $\sigma_t\in]t,T]$ such that $\eta_t\geq L_s$ for each $s\in[t,\sigma_t].$ Let $(y_s,z_s,k_s)$ is a solution of RBSDEs with parameter $(g,t,\eta_t,L)$. By (i), we can easily check
 $$(Y_t,Z_t,K_t)=\left\{
\begin{array}{ll}
    (\eta_t,0,k_t), &  s\in[t,\sigma_t], \\
    (y_s,z_s,k_s), &  s\in[0,t[.
  \end{array}
\right.$$
is a solution of RBSDEs with parameter $(g,\sigma_t,\eta_t,L)$. Thus (ii) holds ture. Now we prove (ii) implies (i). For each $\eta\in {\cal{S}^{\mathrm{2}}}(0,T;{\mathbf{R}})$ satisfying $\eta_t> L_t$ for each $t\in[0,T]$ and each $t\in[0,T[$, we can find a stopping time $\sigma_t$ by setting $$\sigma_t:=\inf\{s\geq t:\eta_t\leq L_s\}\wedge(T-t).$$
Clearly,  $\sigma_t\in]t,T]$ is a stopping time such that $\eta_t>L_s$ for each $s\in[t,\sigma_t[.$ By (ii), there exists a solution $(Y_t,Z_t,K_t)$ of RBSDEs with parameter $(g,\sigma_t,\eta_t,L)$ such that for each $s\in[t,\sigma_t],$ $P-a.s.,\ Y_s=\eta_t.$ Since $\eta_t>L_s$ for $s\in[t,\sigma_t[.$ Thus we have $K_{\sigma_t}-K_{t}=0.$ By Corollary 3.4 and 1 in Remark 3, we can deduce for almost every $t\in[0,T[$, $P-a.s.,\ g(t,\eta_t,0)=0.$ Thus, for each $\zeta\in {\cal{S}^{\mathrm{2}}}(0,T;{\mathbf{R}})$ satisfying $\zeta\geq L,$ and almost every $t\in[0,T[$, we have $P-a.s.,\ g(t,\zeta_t+2^{-n},0)=0.$
By (A2), we obtain (i).\ \ \ $\Box$\\

\textbf{Acknowledgements}

The authors would like to thank the anonymous referees for their careful reading and helpful criticisms.

\end{document}